\documentclass[12pt]{amsart}
\usepackage{amssymb,amscd}

\let\frak\mathfrak

\def\>{\relax\ifmmode\mskip.666667\thinmuskip\relax\else\kern.111111em\fi}
\def\<{\relax\ifmmode\mskip-.333333\thinmuskip\relax\else\kern-.0555556em\fi}
\def\vsk#1>{\vskip#1\baselineskip}
\def\vv#1>{\vadjust{\vsk#1>}\ignorespaces}
\def\vvn#1>{\vadjust{\nobreak\vsk#1>\nobreak}\ignorespaces}

\let\alb\allowbreak

\let\Medskip\medskip
\def\medskip{\par\Medskip}
\let\Bigskip\bigskip
\def\bigskip{\par\Bigskip}

\let\Maketitle\maketitle
\def\maketitle{\Maketitle\thispagestyle{empty}\let\maketitle\empty}

\newtheorem{thm}{Theorem}[section]
\newtheorem{cor}[thm]{Corollary}
\newtheorem{lem}[thm]{Lemma}

\numberwithin{equation}{section}

\theoremstyle{definition}
\newtheorem*{rem}{Remark}

\let\mc\mathcal
\let\nc\newcommand

\nc{\on}{\operatorname}
\nc{\Z}{{\mathbb Z}}
\nc{\C}{{\mathbb C}}
\nc{\N}{{\mathbb N}}
\nc{\pone}{{\mathbb C}{\mathbb P}^1}
\nc{\arr}{\rightarrow}
\nc{\larr}{\longrightarrow}
\nc{\al}{\alpha}
\nc{\W}{{\mc W}}
\nc{\la}{\lambda}
\nc{\su}{\widehat{{\mathfrak sl}}_2}
\nc{\g}{{\mathfrak g}}
\nc{\h}{{\mathfrak h}}
\nc{\m}{{\mathfrak m}}
\nc{\n}{{\mathfrak n}}
\nc{\Gm}{\Gamma}
\nc{\La}{\Lambda}
\nc{\gl}{\widehat{\mathfrak{gl}_2}}
\nc{\bi}{\bibitem}
\nc{\om}{\omega}
\nc{\Res}{\on{Res}}
\nc{\gm}{\gamma}
\nc{\Om}{\Omega}

\def\z{\mathfrak z}

\def\ch{\on{ch}}

\def\End{\on{End}}

\def\Hom{\on{Hom}}
\def\Res{\on{Res}}
\def\rdet{\on{rdet}}
\def\Wr{\on{Wr}}
\def\tbigoplus{\mathop{\textstyle{\bigoplus}}\limits}

\def\B{{\mc B}}
\def\D{{\mc D}}

\def\F{{\mc F}}

\def\Sc{{\mc S}}
\def\V{{\mc V}}

\let\dl\delta
\let\Dl\Delta

\let\si\sigma
\let\Si\Sigma
\let\Sig\varSigma

\let\der\partial

\let\ge\geqslant
\let\geq\geqslant
\let\le\leqslant
\let\leq\leqslant

\nc{\gln}{\mathfrak{gl}_N}
\nc{\sln}{\mathfrak{sl}_N}

\def\glnt{\gln[t]}

\def\Uglnt{U(\glnt)}

\def\slnt{\sln[t]}

\def\beq{\begin{equation}}
\def\eeq{\end{equation}}
\def\be{\begin{equation*}}
\def\ee{\end{equation*}}

\nc{\bean}{\begin{eqnarray}}
\nc{\eean}{\end{eqnarray}}
\nc{\bea}{\begin{eqnarray*}}
\nc{\eea}{\end{eqnarray*}}
\nc{\bs}{\boldsymbol}
\nc{\Ref}[1]{{\rm(\ref{#1})}}

\nc{\Oml}{{\Om_{\bs\la}}}
\nc{\OmLb}{{\Om_{\bs\La,\bs\la,\bs b}}}
\nc{\Ol}{{\mc O_{\bs\la}}}
\nc{\OLb}{{\mc O_{\bs\La,\bs\la,\bs b}}}
\nc{\Ml}{{\mc M_{\bs\la}}}
\nc{\Mlb}{{\mc M_{\bs\La,\bs\la,\bs b}}}
\nc{\Blb}{{\B_{\bs\La,\bs\la,\bs b}}}
\nc{\Omn}{{\Omega_{\bs n,\bs b,\bs K}}}
\nc{\Omlb}{{\bar\Om_{\bs\la}}}
\nc{\VSl}{{(\V^S)_{\bs\la}}}

\nc{\Dlb}{\Dl_{\bs\La,\bs\la,\bs b,\bs K}}
\nc{\ep}{\epsilon}
\nc{\Vn}{{V^{\otimes n}}}
\nc{\Il}{{\mc I_{\bs\la}}}
\nc{\bla}{{\bs\la}}
\nc{\Fla}{\F_{\bs\la}}
\nc{\GL}{{GL_n(\C)}}
\nc{\ga}{\gamma}
\nc{\Ga}{\Gamma}

\begin{document}

\hrule width0pt
\vsk->

\title[Cohomology of a flag variety as a Bethe algebra]
{Cohomology of a flag variety as a Bethe algebra}

\author
[R.\,Rim\'anyi, V.\,Schechtman, V.\,Tarasov, A.\,Varchenko]
{R.\,Rim\'anyi$\>^{\star,1}$, V.\,Schechtman$\>^*$,
V.\,Tarasov$\>^\diamond$, A.\,Varchenko$\>^{\star,2}$}

\maketitle

\begin{center}
{\it $^\star\<$Department of Mathematics, University of North Carolina
at Chapel Hill\\ Chapel Hill, NC 27599-3250, USA\/}
\medskip

{\it $\kern-.4em^*\<$
Institut de Math\'ematiques de Toulouse, 118 Route de Narbonne,
Universit\'e Paul Sabatier, 31062 Toulouse, France}

\medskip
{\it $\kern-.4em^{\diamond}\<$Department of Mathematical Sciences,
Indiana University\,--\>Purdue University Indianapolis\kern-.4em\\
402 North Blackford St, Indianapolis, IN 46202-3216, USA\/}

\medskip
{\it $^\diamond\<$St.\,Petersburg Branch of Steklov Mathematical Institute\\
Fontanka 27, St.\,Petersburg, 191023, Russia\/}
\end{center}

{\let\thefootnote\relax
\footnotetext{\vsk-.8>\noindent
$^1\<${\sl E\>-mail}:\enspace rimanyi@email.unc.edu\>,
supported in part by NSA grant CON:H98230-10-1-0171\\
$^*\<${\sl E\>-mail}:\enspace schechtman@math.ups-tlse.fr\\
$^\diamond\<${\sl E\>-mail}:\enspace vt@math.iupui.edu\>, vt@pdmi.ras.ru\>,
supported in part by NSF grant DMS-0901616\\
$^2\<${\sl E\>-mail}:\enspace anv@email.unc.edu\>,
supported in part by NSF grant DMS-0555327}}

\bigskip
\hfill {\it To the memory of \>V.\,I.\>Arnold}

\bigskip

\medskip
\begin{abstract} We interpret the equivariant cohomology $H^*_{GL_n}(\Fla,\C)$
of a partial flag variety $\Fla$ parametrizing chains of subspaces
$0=F_0\subset F_1\subset\dots\subset F_N = \C^n$,
$\dim F_i/F_{i-1}=\la_i$, as the Bethe algebra $\B^\infty(\V^\pm_\bla)$ of
the $\gln$-weight subspace $\V^\pm_\bla$ of a $\glnt$-module $\V^\pm$.
\end{abstract}

\setcounter{footnote}{0}
\renewcommand{\thefootnote}{\arabic{footnote}}

\section{Introduction}

A Bethe algebra of a quantum integrable model is a commutative algebra of
linear operators (Hamiltonians) acting on the space of states of the model.
An interesting problem is to describe the Bethe algebra as the algebra
of functions on a suitable scheme. Such a description can be considered as
an instance of the geometric Langlands correspondence, see \cite{MTV2},
\cite{MTV3}. The $\gln$ Gaudin model is an example of a quantum integrable
model \cite{G1}, \cite{G2}. The Bethe algebra $\B^K$ of the $\gln$ Gaudin model
is a commutative subalgebra of the current algebra $U(\glnt)$. The algebra
$\B^K$ depends on the parameters $K=(K_1,\dots,K_N)\in \C^N$. Having
a $\glnt$-module $M$, one obtains the commutative subalgebra
$\B^K(M)\subset\End(M)$ as the image of $\B^K$. The geometric interpretation
of the algebra $\B^K(M)$ as the algebra of functions on a scheme leads to
interesting objects. For example, the Bethe algebra
$\B^{K=0}((\otimes_{s=1}^n L_{\bs\La_s}(z_s))^{sing}_\bla)$ of the subspace
of singular vectors of the $\gln$-weight $\bla$ of the tensor product of
finite-dimensional evaluation modules $\otimes_{s=1}^n L_{\bs\La_s}(z_s)$ is
interpreted as the space of functions on the intersection of suitable Schubert
cycles in a Grassmannian variety, see \cite{MTV2}. This interpretation gives
a relation between representation theory and Schubert calculus useful in both
directions.

One of the most interesting $\glnt$-modules is the vector space
$\V = V^{\otimes n}\otimes\C[z_1,\dots,z_n]$ of $V^{\otimes n}$-valued
polynomials in $z_1,\dots,z_n$, where $V=\C^N$ is the standard vector
representation of $\gln$.
The Lie algebra $\glnt$ naturally acts on $\V$ as well as the
symmetric group $S_n$, which permutes the factors of $V^{\otimes n}$ and
variables $z_1,\dots,z_n$ simultaneously. We denote by $\V^+$ and $\V^-$
the $S_n$-invariant and antiinvariant subspaces of $\V$, respectively.
The actions of $\glnt$ and $S_n$ on $\V$ commute, so $\V^+$ and $\V^-$ are
$\glnt$-submodules of $\V$. The Bethe algebra $\B^K$
preserves the $\gln$-weight decompositions $\V^+ = \oplus_\bla \V^+_\bla$ and
$\V^- = \oplus_\bla \V^-_\bla$, $\bla=(\la_1,\dots,\la_N)\in\Z^N_{\geq 0},
|\bla| =n$. The Bethe algebra $\B^K(\V^+_\bla)$ was described in \cite{MTV3}
as the algebra of functions on a suitable space of quasiexponentials
$\{e^{K_iu}(u^{\la_i}+\Sig_{i1}u^{\la_i-1}+\dots+\Sig_{i\la_i})$,
$i=1,\dots,N\}$.
In this paper we give a similar description for $\B^K(\V^-_\bla)$ and study the
limit of the algebras $\B^K(\V^+_\bla)$, $\B^K(\V^-_\bla)$ as all coordinates
of the vector $K$ tend to infinity so that $K_i/K_{i+1}\to\infty$ for all $i$.
We show that in this limit both Bethe algebras $\B^\infty(\V^+_\bla)$,
$\B^\infty(\V^-_\bla)$ can be identified with the algebra of the equivariant
cohomology $H^*_{GL_n}(\Fla,\C)$ of the partial flag variety $\Fla$
parametrizing chains of subspaces
\be
0=F_0\subset F_1\subset\dots\subset F_N = \C^n,
\ee
$\dim F_i/F_{i-1}=\la_i$. This identification was motivated for us by
the considerations in \cite{RV}, \cite{RSV}, where the equivariant cohomology
of the partial flag varieties were used to construct certain conformal blocks
in $V^{\otimes n}$.

Our identification of the Bethe algebra with the algebra of multiplication operators of the equivariant
cohomology $H^*_{GL_n}(\Fla,\C)$
can be considered as a degeneration of the recent description in \cite{O} of the equivariant quantum cohomology
of the partial flag varieties as the Bethe algebra of a suitable Yangian model associated with $V^{\otimes n}$,
cf.~\cite{BMO}.

\medskip
In Section \ref{alg sec} we introduce the Bethe algebra. Section \ref{sec EQ}
contains the main results --- Theorems \ref{main thm}, \ref{lem shap nondeg}.
Theorem \ref{main thm} identifies the algebra of equivariant cohomology
$H^*_{GL_n}(\F_\bla,\C)$ and the Bethe algebras $\B^\infty(\V^+_\bla)$,
$\B^\infty(\V^-_\bla)$. Theorem \ref{lem shap nondeg} says that
the Shapovalov pairing of $\V^+_\bla$ and $\V^-_\bla$ is nondegenerate.
In Section \ref{Spaces of quasi-exponentials} we show that the isomorphisms
of Theorem \ref{main thm} are limiting cases of a geometric Langlands
correspondence. In Section \ref{sec quantum} we explain how the Bethe algebras
$\B^\infty(\V^+_\bla)$, $\B^\infty(\V^-_\bla)$ are related to the quantum
equivariant cohomology $QH_{GL_n\times \C^*}(T^*\F_\bla)$ of the cotangent
bundle $T^*\F_\bla$ of the flag variety $\F_\bla$.
Appendix contains the topological description of $\glnt$-actions on
$\oplus_\bla H^*_{GL_n}(\Fla,\C)$.

\medskip

We thank S.\,Loktev and A.\,Okounkov for useful discussions.

\section{Representations of current algebra $\glnt$}
\label{alg sec}
\subsection{Lie algebra $\gln$}
Let $e_{ij}$, $i,j=1,\dots,N$, be the standard generators of the Lie algebra
$\gln$ satisfying the relations
$[e_{ij},e_{sk}]=\dl_{js}e_{ik}-\dl_{ik}e_{sj}$.
We denote by $\h\subset\gln$ the subalgebra generated by
$e_{ii},\,i=1,\dots,N$. For a Lie algebra $\g\,$, we denote by $U(\g)$ the universal enveloping algebra
of $\g$.

A vector $v$ of a $\gln$-module $M$ has weight
$\bla=(\la_1,\dots,\la_N)\in\C^N$ if $e_{ii}v=\la_iv$ for $i=1,\dots,N$.
We denote by $M_\bla\subset M$ the weight subspace of weight $\bla$.

Let $V=\C^N$ be the standard vector representation of $\gln$ with basis
$v_1,\dots,v_N$ such that $e_{ij}v_k=\delta_{jk}v_i$ for all $i,j,k$.
A tensor power $\Vn$ of the vector representation has a basis
given by the vectors $v_{i_1}\otimes\dots\otimes v_{i_n}$, where
\;$i_j\in\{1,\dots,N\}$. Every such sequence $(i_1,\dots,i_n)$ defines
a decomposition $I=(I_1,\dots,I_N)$ of $\{1,\dots,n\}$ into disjoint subsets
$I_1,\dots,I_N$: \;$I_j=\{k\ |\ i_k=j\}$. We denote the basis vector
$v_{i_1}\otimes\dots\otimes v_{i_n}$ by $v_I$.

Let
\be
\Vn=\!\bigoplus_{\bla\in\Z^N_{\geq 0},\,|\bla|=n}\!(\Vn)_\bla
\ee
be the weight decomposition. Denote $\Il$ the set of all indices $I$ with
$|I_j|=\la_j$, \;$j=1,\dots N$. The vectors $\{v_I, I\in\Il\}$ form a basis of
$(\Vn)_\bla$. The dimension of $(\Vn)_\bla$ equals the multinomial coefficient
$d_\bla:=\frac{n!}{\la_1!\dots\la_N!}$.

\medskip
Let $\Sc$ be the bilinear form on $\Vn$ such that the basis $\{v_I\}$ is
orthonormal. We call $\Sc$ the Shapovalov form.

\subsection{Current algebra $\glnt$}
Let $\glnt=\gln\otimes\C[t]$ be the Lie algebra of $\gln$-valued polynomials
with pointwise commutator.
We identify $\gln$ with the subalgebra $\gln\otimes1$
of constant polynomials in $\glnt$. Hence any $\glnt$-module has the canonical
structure of a $\gln$-module.

The Lie algebra $\glnt$ has a basis $e_{ij}\otimes t^r$, $i,j=1,\dots,N$,
$r\in\Z_{\ge0}$, such that
\vvn.2>
\be
[e_{ij}\otimes t^r,e_{sk}\otimes t^p]\,=\,
\dl_{js}e_{ik}\otimes t^{r+p}-\dl_{ik}e_{sj}\otimes t^{r+p}\,.
\vv.2>
\ee
It is convenient to collect elements of $\glnt$ in generating series
of a variable $u$. For $g\in\gln$, set
$g(u)=\sum_{s=0}^\infty (g\otimes t^s)u^{-s-1}$.

The subalgebra $\z_N[t]\subset\glnt$ with basis
$\sum_{i=1}^Ne_{ii}\otimes t^r,$ $r\in \Z_{\geq 0}$, is central.

\subsection{The $\glnt$-modules $\V^\pm$}
\label{VS}
Let $S_n$ be the permutation group on $n$ elements.
For an $S_n$-module $M$ we denote by
$M^+$ (resp. $M^-$) the subspace of $S_n$-invariants
(resp. antiinvariants).

The group $S_n$ acts on $\C[\bs z]:=\C[z_1,\dots,z_n]$ by permuting
the variables. Denote by $\si_s(\bs z)$, $s=1,\alb\dots,n$,
the $s$th elementary symmetric polynomial in $z_1,\dots,z_n$.

Let $\V$ be the vector space of polynomials in variables $z_1,\dots,z_n$
with coefficients in $V^{\otimes n}$:
\vvn.3>
\be
\V\>=\,V^{\otimes n}\<\otimes_{\C}\C[\bs z]\,.
\vv.2>
\ee
The symmetric group $S_n$ acts on $\V$ by permuting the factors
of $V^{\otimes n}$ and the variables $z_1,\dots,z_n$ simultaneously,
\vvn.2>
\be
\si\bigl(v_1\otimes\dots\otimes v_n\otimes p(z_1,\dots,z_n)\bigr)\,=\,
v_{(\si^{-1})_1}\!\otimes\dots\otimes v_{(\si^{-1})_n}\otimes p(z_{\si_1},\dots,z_{\si_n})\,,\quad\si\in S_n\,.
\kern-3em
\vv.2>
\ee
We are interested in the subspaces $\V^+, \V^-\subset \V$ of $S_n$-invariants
and antiinvariants.

The space $\V$ is a $\glnt$-module,
\vvn-.1>
\be
g\otimes t^r\,\bigl(v_1\otimes\dots\otimes v_n\otimes p(\bs z))\,=\,
\sum_{s=1}^n
v_1\otimes\dots\otimes gv_s\otimes\dots\otimes v_n\otimes z_s^rp(\bs z)\, .
\vv.2>
\ee
The image of the subalgebra $U(\z_N[t])\subset \Uglnt$ in
$\End(\V)$ is the algebra of operators of multiplication
by elements of \,$\C[\bs z]^+$.
The $\glnt$-action on $\V$ commutes with the $S_n$-action.
Hence $\V^+$ and $\V^-$ are $\glnt$-submodules of $\V$.
The subspaces $\V^+$ and $\V^-$ are free $\C[\bs z]^+$-modules of rank $N^n$.

Consider the $\gln$-weight decompositions
\vvn.3>
\be
\V^+=\oplus_{\bla\in\Z^N_{\geq 0}, |\bla|=n} \V_\bla^+,
\qquad \V^-=\oplus_{\bla\in\Z^N_{\geq 0}, |\bla|=n} \V_\bla^- .
\vv.3>
\ee
For any $\bla$, the subspaces $\V_\bla^+$ and $\V_\bla^-$ are free
$\C[\bs z]^+$-modules of rank $d_\bla$.
\medskip

Denote by $\frac 1D\V^-$ the vector space of all $V^{\otimes n}$-valued
\vvn.16>
rational functions of the form
$\frac 1D x$, $x\in\V^-$, $D=\prod_{1\leq i<j\leq n}(z_j-z_i)$.
The Shapovalov form induces a $\C[\bs z]^+$-bilinear map
\vvn.2>
\be
\Sc_{+-} : \V^+\otimes \frac 1D\V^-\to\,\C[\bs z]^+ .
\vv.3>
\ee
The $\glnt$-module structures on $\V^+$ and $\frac 1D\V^-$ are contravariantly related through the Shapovalov form,
\vvn-.3>
\be
\Sc_{+-}\bigl((e_{ij}\otimes t^r)x,\frac 1D\,y\bigr)\,=\,
\Sc_{+-}\bigl(x,(e_{ji}\otimes t^r)\frac 1D\,y\bigr)
\qquad\text{for all \;} i, j, x, y.
\ee

\subsection{Bethe algebra}
\label{bethesec}

Given an $N\times N$ matrix $A$ with possibly noncommuting entries $a_{ij}$,
we define its row determinant to be
\vvn.1>
\be
\rdet A\,=
\sum_{\;\si\in S_N\!} (-1)^\si\,a_{1\si(1)}a_{2\si(2)}\dots a_{N\si(N)}\,.
\vv-.2>
\ee

Let $K=(K_1,\dots,K_N)$ be a sequence of complex numbers.
Let $\der$ be the operator of differentiation in a variable $u$.
Define the universal differential operator \,$\D^K$ by
\vvn.3>
\be
\D^K=\,\rdet\left( \begin{matrix}
\der-K_1-\>e_{11}(u) & -\>e_{21}(u)& \dots & -\>e_{N1}(u)\\
-\>e_{12}(u) &\der-K_{2}-e_{22}(u)& \dots & -\>e_{N2}(u)\\
\dots & \dots &\dots &\dots \\
-\>e_{1N}(u) & -\>e_{2N}(u)& \dots & \der-K_{N}-e_{NN}(u)
\end{matrix}\right).
\vv.1>
\ee
It is a differential operator in the variable $u$, whose coefficients are
formal power series in $u^{-1}$ with coefficients in
$\Uglnt$,
\vvn.1>
\be
\D^K=\,\der^N+\sum_{i=1}^N\,B^K_i(u)\,\der^{N-i}\>,
\qquad
B_i^K(u)\,=\,\sum_{j=0}^\infty B^K_{ij}\>u^{-j}
\vv-.1>
\ee
and \,$B^K_{ij}\in\Uglnt$ \,for \,$i=1,\dots,N$, \,$j\geq 0$.

Denote by $\B^K$ the unital subalgebra of $\Uglnt$ generated by $B^K_{ij}$
with $i=1,\dots,N$, $j\geq 0$. The subalgebra $\B^K$ is called the Bethe
algebra with parameters $K$.

\begin{thm}
%[\cite{T}, \cite{CT}, \cite{MTV1}]
\label{T-thm}
The algebra $\B^K$ is commutative. The algebra $\B^K$ commutes with
the subalgebra $U(\h)\subset \Uglnt$. If $K=0$, then the algebra $\B^{K=0}$
commutes with the subalgebra $U(\gln)\subset \Uglnt$.
\qed
\end{thm}

The theorem is proved for $K=0$ in \cite{T} and for nonzero $K$ in \cite{CT},
\cite{MTV1}.

Each element $B^K_{ij}$ is a polynomial in $K_1,\dots,K_N$. We define
$\B^\infty$ to be the unital subalgebra of $\Uglnt$ generated by the leading
terms of the elements $B^K_{ij}$, $i=1,\dots,N$, \,$j\geq 0$, as $K$ tends to
infinity so that $K_i/K_{i+1}\to \infty$ for all $i$.

\begin{lem}
\label{lem infty}
The algebra \;$\B^\infty$ is the unital subalgebra generated by the elements
$e_{ii}\otimes t^j$ with $i=1,\dots,N$, \;$j\geq 0$.
\end{lem}
\begin{proof}
We have \;$B_{i0}^K\,=\,(-1)^i\,K_1\dots K_i\,\bigl(\<\>1+o(1)\bigr)$\,,
\;and
\vvn.3>
\be
B^K_{ij}\,=\,(-1)^i\,K_1\dots K_{i-1}\,
\Bigl(\,\sum_{m=i}^N e_{mm}\otimes t^{j-1}+o(1)\Bigr)
\ee
for \,$j>0$\,, where \,$o(1)$ stands for the terms vanishing as \,$K$ tends
to infinity.
\end{proof}

\begin{rem}
There are $N\<\>!$ asymptotic zones labeled by elements of $S_N$ in which $K$
may tend to infinity. For $\si\in S_N$ we may assume that all coordinates of
$K$ tend to infinity and $K_{\si_i}/K_{\si_{i+1}} \to \infty$ for all $i$.
It is easy to see that the limiting Bethe algebra $\B^\infty$ does not depend
on $\si$.
\end{rem}

The algebra $\B^\infty$ is commutative and contains $U(\z_N[t])$.
The algebra $\B^\infty$ commutes with the subalgebra $U(\h)\subset \Uglnt$.

\medskip
As a subalgebra of $\Uglnt$, the Bethe algebra $\B^K$ acts on any
$\glnt$-module $M$. Since $\B^K$ commutes with $U(\h)$, it preserves the weight
subspaces $M_\bla$. If $K=0$, then $\B^{K=0}$ preserves the singular weight
subspaces $M_\bla^{sing}$. We will study the action of $\B^\infty$ on
the weight subspaces $\V^+_\bla$, $\V^-_\bla$.

\begin{lem}
\label{eii}
The element \,$\sum_{i=1}^N e_{ii}\otimes t^r \in U(\z_N[t])$ acts on \,$\V$
as the operator of multiplica\-tion by \,$\sum_{s=1}^n z_s^r$.
\qed
\end{lem}

If $L\subset M$ is a $\B^K$-invariant subspace, then the image of $\B^K$ in
$\End(L)$ will be called the Bethe algebra of $L$ and denoted by $\B^K(L)$.

\section{Equivariant cohomology of partial flag varieties}
\label{sec EQ}

\subsection{Partial flag varieties}
\label{sec Partial flag varieties}
For $\bla\in\Z^N_{\geq 0}, \;|\bla|=n$, consider the partial flag variety
$\Fla$ parametrizing chains of subspaces
\be
0=F_0\subset F_1\subset \dots\subset F_N = \C^n
\ee
with $\dim F_i/F_{i-1}=\la_i$, \,$i=1,\dots,N$.

Let $T^n\subset GL_n$ be the torus of diagonal matrices.
The group $T^n\subset GL_n(\C)$ acts on $\C^n$ and hence on $\Fla$.
The set of fixed points $\Fla^{T^n}$ of the torus action consists of
coordinate flags
$F_I = (F_0\subset\dots\subset F_N)$, \,$I=(I_1,\dots,I_N) \in \Il$,
where $F_i$ is the span of the basis vectors $v_j\in\C^n$
with $j\in I_1\cup\dots\cup I_i$.
The fixed points are in a one-to-one
correspondence with the elements of $\Il$ and hence with the basis vectors
of $V_\bla$.

We consider the $\GL$-equivariant cohomology
\be
H_\bla\,=\,H^*_{GL_n}(\Fla,\C)\,.
\ee
Denote by $\Gamma_i=\{\ga_{i1},\dots,\ga_{i\la_i}\}$ the set of the Chern roots of the bundle over $\Fla$
with fiber $F_i/F_{i-1}$. Denote by $\bs z=\{z_1,\dots,z_n\}$ the Chern roots corresponding to the factors
of the torus $T^n$. Then
\beq
\label{relations}
H_\bla\,=\,\C[\bs z; \Ga_1;\dots;\Ga_N]^{S_n\times S_{\la_1}\times\dots\times
S_{\la_N}}\Big/\Bigl\langle\,
\prod_{i=1}^N\prod_{j=1}^{\la_i}\,(1+u\ga_{ij})\,=\,\prod_{i=1}^n\,(1+uz_i)
\Bigr\rangle\,.
\eeq
The cohomology $H_\bla$ is a module over $H^*_{GL_n}({pt},\C) = \C[\bs z]^+$.

Let $J_H\subset H_\bla$ be the ideal generated by the polynomials
$\si_i(\bs z)$, \,$i=1,\dots,n$. Then \,$H_\bla/J_H=H^*(\Fla,\C)$\,.

\subsection{Integration over $\Fla$}
We will need the integration map $\int : H_\bla\to H^*_{GL_n}({pt},\C)$.
The following formula~\Ref{int formula} gives the integration map in terms of
the fixed point set $\Fla^{T^n}$.

For a subset $A\subset \{1,\dots,N\}$ denote $\bs z_A=\{z_a, \,a\in A\}$.
For $I=(I_1,\dots,I_N)\in\Il$ denote
\be
R(\bs z_{I_1}|\bs z_{I_2}|\ldots |\bs z_{I_m})\,=\,
\prod_{i<j}\,\prod_{a\in I_i,\,b\in I_j}\!(z_b - z_a).
\ee
The Atiyah-Bott equivariant localization theorem \cite{AB} says that
for any $[h(\bs z, \Ga_1,\dots,\Ga_N)] \in H_\bla$,
\vvn-.2>
\beq
\label{int formula}
\int [h]\,=\,\sum_{I\in\Il} \frac {h(\bs z, \bs z_{I_1},\dots,\bs z_{I_N})}
{R(\bs z_{I_1}|\bs z_{I_2}|\ldots |\bs z_{I_N})}\;.
\vv.2>
\eeq
More precisely, it is enough to verify this identity after the base change
$H_{GL_n}^*(\F_\bla) \to H_{T^n}^*(\F_\bla)$, where $T^n\subset GL_n$ is
the maximal torus, and in the $T^n$-equivariant cohomology this is
the Atiyah-Bott localization formula. Note that the expression for $R$ is
nothing else but
the Euler class of the tangent space at a $T^n$-fixed point, see formula (2.5)
in \cite{RSV}.

Clearly, the right hand side in \Ref{int formula} lies in $\C[\bs z]^+$.
The integration map induces the pairing
\be
(\,,)\,:\,H_\bla\otimes H_\bla\to\C[\bs z]^+,\qquad
[h]\otimes [g] \mapsto \int [hg]\;.
\ee
After factorization by the ideal $J_H$ we obtain the nondegenerate Poincare
pairing
\vvn.3>
\be
(\,,)\,:\,H^*(\Fla,\C)\otimes H^*(\Fla,\C) \ \to\ \C .
\ee

\subsection{$H_\bla$ and $\V^\pm$}

\begin{lem}
\label{i lem}
The maps
\begin{align*}
i^+_\bla\,:\,H_\bla \to \V^+_\bla, &\qquad
[h(\bs z, \Ga_1,\dots,\Ga_N)]\,\mapsto\,
\sum_{I\in\Il} v_I\otimes h(\bs z, \bs z_{I_1},\dots,\bs z_{I_N}) ,
\\
i^-_\bla\,:\,H_\bla \to \frac 1D\V^-_\bla, &\qquad
[h(\bs z, \Ga_1,\dots,\Ga_N)]\,\mapsto\,
\sum_{I\in\Il} v_I\otimes \frac
{h(\bs z, \bs z_{I_1},\dots,\bs z_{I_N})}
{R(\bs z_{I_1}|\bs z_{I_2}|\ldots |\bs z_{I_N})}
\notag
\end{align*}
are well-defined isomorphisms of $\C[\bs z]^+$-modules.
\qed
\end{lem}

\begin{proof}
If $h$ belongs to the ideal of relations in \Ref{relations} then
$h(z,z_{I_1},\dots,z_{I_N})=0$ for any $I$, because the $\Gamma_i=z_{I_i}$
substitution makes the generators of the ideal identities.
This proves well-definedness.

Consider the $\C[\bs z]^+$-module
$\C[\bs z]^{S_{\la_1}\times \ldots \times S_{\la_N}}$ of polynomials symmetric
in the first $\lambda_1$ variables, the next $\lambda_2$ variables, etc.
In Schubert calculus it is known that this module is free of rank $d_\bla$,
and that it is isomorphic to $H_\bla$ under the correspondence
\beq
\label{rr_isom}
p\in \C[\bs z]^{S_{\la_1}\times \ldots \times S_{\la_N}}\,\,\leftarrow\!\!
\!\!\to\,\,\bigl[\>p(\Gamma_1,\ldots,\Gamma_N)\bigr]\in H_\bla.
\eeq
An element $\sum_{I\in\Il} v_I \otimes p_I(\bs z)$ of $\V_\bla$ belongs to $\V^+_\bla$, if and only if
$p_I(\bs z)=p(z_{I_1},\ldots,z_{I_N})$ for a polynomial $p\in \C[\bs z]^{S_{\la_1}\times \ldots \times S_{\la_N}}$. This shows that $\V^+_\bla$ is isomorphic to $\C[\bs z]^{S_{\la_1}\times \ldots \times S_{\la_N}}$, and that $i^+_\bla$ is the composition of this isomorphism with (\ref{rr_isom}).

A similar argument shows that $i^-_\bla$ is also an isomorphism.
\end{proof}

\begin{cor}
\label{lem Shap Poinc}
The Shapovalov form and the Poincare pairing are related by the formula
\be
\Sc_{+-}(i_+[h], i_-[g]) = \int [h][g]\;.
\vv->
\ee
\qed
\end{cor}

Let $A$ be a commutative algebra. The algebra $A$ considered as an $A$-module
is called the regular representation of $A$. Here is our main result.

\begin{thm}
\label{main thm}
\strut
\begin{enumerate}
\item[(i)]
The maps \,\,$\xi^\pm_\bla:e_{ii}\otimes t^r|_{\V^\pm_\bla}\mapsto\,
\sum_{j=1}^{\la_i} \gamma_{ij}^r$ \,define isomorphisms of the algebras
\,$\B^\infty(\V^\pm_\bla)$ and \,$H_\bla$\,.
\item[(ii)]
The maps \;$\xi^+_\bla,\,i^+_\bla$ identify the \,$\B^\infty(\V^+_\bla)$-module
$\V^+_\bla$ with the regular representation of $H_\bla$.
\item[(iii)]
The maps \;$\xi^-_\bla,\,i^-_\bla$ identify the \,$\B^\infty(\frac 1D\V^-_\bla)$-module
$\frac 1D\V^-_\bla$ with the regular representation of $H_\bla$.
\end{enumerate}
\end{thm}

The theorem follows from Lemmas \ref{lem infty}, \ref{eii} and \ref{i lem}.

\subsection{Cohomology as $\glnt$-modules}
\label{rr_sec}
Let $J$ be the ideal of $\C[\bs z]^+$ generated by the ele\-mentary symmetric
functions $\si_i(\bs z)$, \,$i=1,\dots,n$. Define $J^+=J\V^+$ and
$J^-=\frac 1DJ\V^-$. Clearly, $J^+$ is a $\glnt$-submodule of $\V^+$ and $J^-$
is a $\glnt$-submodule of $\frac 1D\V^-$. The $\glnt$-module $\V^+/J^+$ is
graded and has dimension $N^n$ over $\C$, see \cite{MTV2}. Similarly,
$\frac 1D\V^-/J^-$ is a graded $\glnt$-module of the same dimension.

\begin{thm}
\label{lem shap nondeg}
The Shapovalov form establishes a nondegenerate pairing
\vvn.2>
\be
\Sc_{+-} \ :\ \V^+/J^+ \otimes \frac 1D \V^-/J^- \ \to \ \C .
\vv.3>
\ee
\end{thm}

\noindent
The theorem follows from Lemmas \ref{i lem}, \ref{lem Shap Poinc} and
the nondegeneracy of the Poincare pairing.

\begin{cor}
\label{cor Shap rela}
The $\glnt$-modules $\V^+/J^+$ and $ \frac 1D \V^-/J^-$ are contravariantly
related through the Shapovalov form,
\,$\Sc_{+-}\bigl((e_{ij}\otimes t^r)x,\frac 1D y\bigr)
= \Sc_{+-}\bigl(x,(e_{ji}\otimes t^r)\frac 1D y\bigr)$ \,for all \,$i,j,x,y$.
\end{cor}

Let $W_n$ be the $\glnt$-module generated by a vector $w_n$ with the
defining relations:
\vvn-.2>
\begin{alignat*}2
e_{ii} &(u)w_n=\>\dl_{1i}\,\frac nu\,w_n\,, && i=1,\dots,N\,,
\\[4pt]
e_{ij} & (u)w_n=\>0\,, && 1\le i<j\le N\,,
\\[4pt]
(e_{ji} &{}\otimes1)^{n\dl_{1i}+1}w_n=\>0\,,\qquad && 1\le i<j\le N\,.
\\[-12pt]
\end{alignat*}
As an $\slnt$-module, the module $W_n$ is isomorphic to the Weyl module from
\cite{CL}, \cite{CP}, corresponding to the weight $n\om_1$, where $\om_1$ is
the first fundamental weight of $\sln$.

In \cite{MTV2} an isomorphism of $\V^+/J^+$ and the Weyl module $W_n$ is
constructed.

\begin{cor}
\label{lem Weyl dual}
The Shapovalov form $\Sc_{+-}$ establishes an isomorphism of
$ \frac 1D \V^-/J^-$ and the contravariantly dual of the Weyl module $W_n$.
\end{cor}

Here is an application of this fact.
For $\bla\in\Z^N_{\geq 0}$, $|\bla|=n$, $\la_1\geq\dots\geq\la_N$, denote
\be
(\frac 1D \V^-/J^-)^{sing}_\bla\,=\,\{v\in \frac 1D \V^-/J^-\ | \ e_{ij}v=0 \
\text{for}\ i<j,\ e_{ii}v=\la_iv\ \;\text{for}\ \,i=1,\dots,N\}\,.
\vv.2>
\ee
This is a graded space. Denote by $\bigl((\frac1D\V^-/J^-)^{sing}_\bla\bigr)_k$
the subspace of all elements of $\bs z$-degree $k$. Define the graded character
by the formula
\be
\ch\bigl((\frac 1D \V^-/J^-)^{sing}_\bla\bigr)\,=\,
\sum_k q^k \dim ((\frac 1D \V^-/J^-)^{sing}_\bla)_k .
\ee

\begin{cor}
\label{cor grading}
We have
\beq
\label{grad}
\ch\bigl((\frac 1D \V^-/J^-)^{sing}_\bla\bigr)\,=\,
\frac{(q)_n\prod_{1\le i<j\le N}(1-q^{\la_i-\la_j+j-i})}
{\prod_{i=1}^N(q)_{\la_i+N-i}}\ q^{\,-\sum_{1\le i<j\le N}\la_i\la_j},
\eeq
where $\,(q)_a=\prod_{j=1}^a(1-q^j)\,$.
\end{cor}

The corollary follows from Lemma 2.2 in \cite{MTV2} and
Corollary \ref{lem Weyl dual}.

\medskip
The isomorphisms
\beq
\label{two isom}
i^+=\,\tbigoplus_\bla\,i^+_\bla \,:\,\tbigoplus_\bla\,H_\bla\,\to\;\V^+,
\qquad
i^-=\,\tbigoplus_\bla\,i^-_\bla\,:\,\tbigoplus_\bla\,H_\bla\,\to\;\frac 1D\V^-
\vv.3>
\eeq
induce two graded $\glnt$-module structures on $\oplus_\bla\,H_\bla$ denoted
by $\rho^+$ and $\rho^-$, respectively. These module structures descend to two
graded $\glnt$-module structures on the cohomology with constant coefficients
\vvn.2>
\be
H(\C)\,:=\!
\bigoplus_{\bla\in\Z^N_{\geq 0},\,|\bla|=n}\!H^*(\Fla,\C)\,,
\vv-.2>
\ee
denoted by the same letters $\rho^+$ and $\rho^-$.

\begin{cor}
\label{cor H(C) modules}
The $\glnt$-module $H(\C)$ with the $\rho^+\!$-structure is isomorphic to the
Weyl module $W_n$. The $\glnt$-module $H(\C)$ with the $\rho^-\!$-structure is
isomorphic to the contravariant dual of the Weyl module $W_n$.
\end{cor}

The $\rho^\pm$ structures can be defined topologically, see \cite{RSV} and
Appendix. The $\rho^-$-structure appears to be more preferable. It was used
in \cite{RV}, \cite{RSV} to construct conformal blocks in the tensor power
$V^{\otimes n}$.

\section{Isomorphisms $i^\pm_\bla$ as a geometric Langlands correspondence}
\label{Spaces of quasi-exponentials}

\subsection{The $\V^+_\bla$ case}
\label{The V+_bla case}
The following geometric description of the $\B^K$-action on $\V^+_\bla$ was
given in \cite{MTV3} as an example of the geometric Langlands correspondence.

Let $K=(K_1,\dots,K_N)$ be a sequence of distinct complex numbers.
Let $\bla\in\Z^N_{\geq 0}$\,, \,$|\bla|=n$.
Introduce the polynomial algebras
\vvn.3>
\be
\C[\bs\Sig]\,:=\,\C[\Sig_{ij},\,i=1,\dots,N,\,j=1,\dots,\la_i]\,,\qquad
\C[\bs\si]\,:=\,\C[\si_1,\dots,\si_n]\,.
\ee
Define
\be
\Sig_i(u)\,=\,e^{K_iu}\,(u^{\la_i}+\Sig_{i1}u^{\la_i-1}+\dots+\Sig_{i\la_i})\,,
\qquad i=1,\dots,N\,.
\ee
For arbitrary functions $g_1(u),\dots,g_N(u)$, introduce the Wronskian
determinant by the formula
\vvn-.3>
\be
\Wr(g_1(u),\dots,g_N(u))\,=\,
\det\left(\begin{matrix} g_1(u) & g_1'(u) &\dots & g_1^{(N-1)}(u) \\
g_2(u) & g_2'(u) &\dots & g_2^{(N-1)}(u) \\ \dots & \dots &\dots & \dots \\
g_N(u) & g_N'(u) &\dots & g_N^{(N-1)}(u)
\end{matrix}\right).
\ee
We have
\be
\Wr(\Sig_1(u),\dots,\Sig_N(u))\,=\,
e^{\>\sum_{i=1}^N K_iu}\!\prod_{1\le i<j\le N}\!(K_j-K_i)\cdot
\Bigl(u^n+\sum_{s=1}^n\,(-1)^s\>A_s^K(\bs\Sig)\,u^{n-s}\Bigr)\,,
\ee
where $A_1^K(\bs\Sig),\dots,A_n^K(\bs\Sig)\in\C[\bs\Sig]$.
Define an algebra homomorphism
\be
\W^K :\,\C[\bs\si] \to \C[\bs\Sig] ,\qquad
\si_s \mapsto A_s^K(\bs\Sig)\,.
\ee
The homomorphism defines a $\C[\bs\si]$-module structure on $\C[\bs\Sig]$.

\medskip
Define a differential operator $\D^K_{\bs\Si}$ by
\vvn-.7>
\be
\D^K_{\bs\Sig}=\,\frac1{\Wr(\Sig_1(u),\dots,\Sig_N(u))}\,\rdet
\left(\begin{matrix} \Sig_1(u) & \Sig_1'(u) &\dots & \Sig_1^{(N)}(u) \\
\Sig_2(u) & \Sig_2'(u) &\dots & \Sig_2^{(N)}(u) \\ \dots & \dots &\dots & \dots \\
1 & \der &\dots & \der^N
\end{matrix}\right).
\ee
It is a differential operator in the variable $u$, whose coefficients are
formal power series in $u^{-1}$ with coefficients in $\C[\bs\Sig]$,
\vvn-.4>
\beq
\label{DOK}
\D^K_{\bs\Sig}=\,\der^N+\sum_{i=1}^N\,F_i^K(u)\,\der^{N-i}\>,
\qquad F_i^K(u)\,=\,\sum_{j=0}^\infty F^K_{ij}\>u^{-j}\,,
\eeq
and $F_{ij}^K\in\C[\bs\Sig]$, \,$i=1,\dots,N$, \,$j\geq 0$.

\begin{thm}[\cite{MTV3}]
\label{thm 1 mtv7}
The map
\be
\tau_\bla^{K+}\,:\,B^K_{ij}|_{\V^+_\bla}\,\mapsto\,F^K_{ij}
\ee
defines an isomorphism of the Bethe algebra
$\B^K(\V^+_\bla)$ and the algebra $\C[\bs\Sig]$.
The isomorphism $\tau_\bla^{K+}$ becomes an isomorphism of the $U(\z_N[t])|_{\V^+_\bla}$-module $\B^K(\V^+_\bla)$
and the $\C[\bs\si]$-module $\C[\bs\Sig]$ if we identify the algebras
$U(\z_N[t])|_{\V^+_\bla}$ and $\C[\bs\si]$ by the map
$\si_s[\bs z]\mapsto \si_s$, \,$s=1,\dots,n$.
\end{thm}

Denote
\vvn-.4>
\be
v^+ =\,\sum_{I\in\Il} v_I \in \V^+_\bla .
\ee

\begin{thm}[\cite{MTV3}]
\label{thm 2 mtv7}
The map
\be
\mu_\bla^{K+}\,:\,
B^K_{ij}v^+\,\mapsto\,F_{ij}^K\,,
\ee
defines a linear isomorphism $\V^+_\bla \to \C[\bs\Sig]$.
The maps $\tau_\bla^{K+}, \mu_\bla^{K+}$ give an isomorphism of
the $\B^K(\V^+_\bla)$-module $\V^+_\bla$ and the regular representation
of the algebra $\C[\bs\Sig]$.
\end{thm}

\subsection{The limit of $\tau_\bla^{K+}$ and $\mu_\bla^{K+}$ as $K\to\infty$}
\label{sec lim K+ infty}
Let all the coordinates of the vector $K$ tend to infinity so that
$K_i/K_{i+1}\to\infty$ for $i=1,\dots,N-1$. Then the homomorphism $\W^K$ has
a limit $\W^\infty$. Namely, define $A^\infty_s(\bs\Sig)$ by
the formula
\be
\prod_{i=1}^N\,(u^{\la_i}+\Sig_{i1}u^{\la_i-1}+\dots+\Sig_{i\la_i})\,=
u^n+\sum_{s=1}^n\,(-1)^s\>A_s^\infty(\bs\Sig)\,u^{n-s}\,.
\ee
Then
\vvn-.4>
\beq
\label{W inft}
\W^\infty :\,\C[\bs\si]\,\to\,\C[\bs \Sig],\qquad
\si_s\,\mapsto\,A_s^\infty(\bs\Sig)\,.
\eeq
Define algebra homomorphisms
\beq
\label{alg isos}
\C[\bs\si]\,\to\,\C[\bs z]^+,\qquad
\eta:\C[\bs\Sig]\,\to\,H_\bla\,,
\eeq
by the agreement that the first one sends $\si_s$ to $\si_s(\bs z)$ for all
$s$\,, and the second one sends $(-1)^s\Sig_{is}$ to the $s$-th elementary
symmetric function of \,$\gm_{i1},\dots,\gm_{i\la_i}$ for all \,$i,s$.
Clearly, the defined maps are isomorphisms.

\begin{lem}
\label{lem iso C[Sigma] H-bla}
The isomorphisms \Ref{alg isos} identify the $\C[\bs\si]$-module
$\C[\bs\Sig]$ defined by formula \Ref{W inft} and the $\C[\bs z]^+$-module
$H_\bla$.
\end{lem}

Let \;$p_i(u)\,=\,u^{\la_i}+\Sig_{i1}u^{\la_i-1}+\dots+\Sig_{i\la_i}$
\,for all \,$i=1,\dots N$. Notice that
\beq
\label{etap}
\eta\bigl(p_i(u)\bigr)\,=\,\prod_{j=1}^{\la_i}\,(u-\gm_{ij})\,,\qquad
\eta\Bigl(\>\frac{p_i'(u)}{p_i(u)}\Bigr)\,=\,
\sum_{r=0}^\infty\,\sum_{j=1}^{\la_i}\,\gm_{ij}^r\,u^{-r-1}\,.
\vv.2>
\eeq

\begin{lem}
\label{Fkij}
We have \;$F_{i0}^K\,=\,(-1)^i\,K_1\dots K_i\,\bigl(\<\>1+o(1)\bigr)$\,, \;and
\vvn.1>
\be
\sum_{j=1}^\infty\,F^K_{ij}u^{-j}\,=\,(-1)^i\,K_1\dots K_{i-1}\,
\Bigl(\,\sum_{m=i}^N\,\frac{p_m'(u)}{p_m(u)}+o(1)\Bigr)\,,
\vv.2>
\ee
where \,$o(1)$ stands for the terms vanishing as \,$K$ tends to infinity.
\end{lem}
\begin{proof}
Let \;$y_i(u)=\Wr\bigl(\Sig_i(u),\dots\Sig_N(u)\bigr)$\,, \;$i=1,\dots N$.
Then the operator \,$\D^K_{\bs\Sig}$ can be factorized:
\vvn.2>
\beq
\label{DKfact}
\D^K_{\bs\Sig}\,=\,
\Bigl(\der-\frac{y_1'(u)}{y_1(u)}+\frac{y_2'(u)}{y_2(u)}\,\Bigr)\,\dots\,
\Bigl(\der-\frac{y_{N-1}'(u)}{y_{N-1}(u)}+\frac{y_N'(u)}{y_N(u)}\,\Bigr)\,
\Bigl(\der-\frac{y_N'(u)}{y_N(u)}\,\Bigr)\,,
\vv.4>
\eeq
see \cite{MV}. Since
\be
y_i(u)\,=\,(-1)^{(N-i)(N-i-1)/2}\,
K_i^{N-i}\dots K_{N-1}\,\bigl(\<\>p_i(u)\dots p_N(u)+o(1)\bigr)
\,e^{\>\sum_{m=i}^NK_mu}
\vv.2>
\ee
as \,$K$ tends to infinity,
the claim follows from formulae~\Ref{DOK} and~\Ref{DKfact}.
\end{proof}

\begin{thm}
\label{thm + K infty}
\strut
\begin{enumerate}
\item[(i)]
The map \;$\eta\circ\<\tau_\bla^{K+}\!:\B^K(\V^+_\bla)\to H_\bla$ \,tends to
\vvn.1>
the isomorphism \;$\xi^+_\bla:\B^\infty(\V^+_\bla)\to H_\bla$, see
Theorem \ref{main thm}, \,as \,$K$ tends to infinity.
\vsk.2>
\item[(ii)]
The map \;$\eta\circ\<\mu_\bla^{K+}\!:\V^+_\bla\to H_\bla$ \,tends to
\vvn.1>
the isomorphism \;$(i^+_\bla)^{-1}:\V^+_\bla\to H_\bla$,
see Lemma \ref{i lem}, \,as \,$K$ tends to infinity.
\end{enumerate}
\end{thm}

\begin{proof}
The statement follows from the definitions of the maps, Lemma~\ref{Fkij},
formulae \Ref{etap}, and the proof of Lemma~\ref{lem infty}.
\end{proof}

\subsection{The $\frac1D\V^-_\bla$ case}
\label{The V-_bla case}

Theorem \ref{lem shap nondeg} allows us to establish a geometric description of
the $\B^K$-action on $\frac 1D\V^-$ which is analogous to the description of
the $\B^K$-action on $\V^+$.

\begin{thm}
\label{thm 1 }
The map
\be
\tau_\bla^{K-}\,:\,B^K_{ij}|_{\frac 1D\V^-_\bla}\,\mapsto\,F^K_{ij}
\vv.3>
\ee
defines an isomorphism of the Bethe algebra
$\B^K(\frac 1D\V^-_\bla)$ and the algebra $\C[\bs\Sig]$.
The isomorphism $\tau_\bla^{K-}$ becomes an isomorphism of the $U(\z_N[t])|_{\frac 1D\V^-_\bla}$-module $\B^K(\frac 1D\V^-_\bla)$
and the $\C[\bs\si]$-module $\C[\bs\Sig]$ if we identify the algebras
$U(\z_N[t])|_{\frac 1D\V^-_\bla}$ and $\C[\bs\si]$ by the map
$\si_s[\bs z]\mapsto \si_s$, \,$s=1,\dots,n$.
\end{thm}

Denote
\vv-.6>
\be
v^-=\,\sum_{I\in\Il} v_I\otimes
\frac1{R(\bs z_{I_1}|\bs z_{I_2}|\ldots |\bs z_{I_N})}\in\frac1D\V^-_\bla\,.
\ee

\begin{thm}
\label{thm 2}
The map
\vv-.2>
\be
\mu_\bla^{K-}\,:\,B^K_{ij}v^-\,\mapsto\,F_{ij}^K,
\vv.3>
\ee
defines a linear isomorphism $\frac 1D\V^-_\bla \to \C[\bs\Sig]$.
The maps $\tau_\bla^{K-}, \mu_\bla^{K-}$ give an isomorphism of
the $\B^K(\frac 1D\V^-_\bla)$-module $\frac 1D\V^-_\bla$ and the regular
representation of the algebra $\C[\bs\Sig]$.
\end{thm}

The proofs of Theorems \ref{thm 1 } and \ref{thm 2} are basically word by word
the same as the proofs of Theorems \ref{thm 1 mtv7} and \ref{thm 2 mtv7} in
\cite{MTV3}.

It is interesting to note that the element $v^-$ becomes a conformal block
under certain conditions and satisfies a {\sl KZ\/} equation with respect
to $\bs z$, see \cite{V}, \cite{RV}, \cite{RSV}.

\subsection{The limit of $\tau_\bla^{K-}$ and $\mu_\bla^{K-}$ as $K\to\infty$}
\label{sec lim K- infty}
Let all the coordinates of the vector $K$ tend to infinity so that
$K_i/K_{i+1}\to \infty$ for $i=1,\dots,N-1$.

\begin{thm}
%\label{thm - K infty}
\strut
\begin{enumerate}
\item[(i)]
The map \;$\eta\circ\<\tau_\bla^{K-}\!:\B^K(\V^-_\bla)\to H_\bla$ \,tends to
\vvn.1>
the isomorphism \;$\xi^-_\bla:\B^\infty(\V^-_\bla)\to H_\bla$, see
Theorem \ref{main thm}, \,as \,$K$ tends to infinity.
\vsk.2>
\item[(ii)]
The map \;$\eta\circ\<\mu_\bla^{K-}\!:\V^-_\bla\to H_\bla$ \,tends to
\vvn.1>
the isomorphism \;$(i^-_\bla)^{-1}:\V^-_\bla\to H_\bla$,
see Lemma \ref{i lem}, \,as \,$K$ tends to infinity.
\end{enumerate}
\end{thm}

The proof is similar to the proof of Theorem \ref{thm + K infty}.

\subsection{The $(\frac 1D \V^-)^{sing}_\bla$ case}
\label{sec sing case}
Formula \Ref{grad} for the graded character of $(\frac 1D\V^-/J^-)^{sing}_\bla$ is the analog of the formula
for the graded character of $(\V^+/J^+)^{sing}_\bla$ in \cite{MTV2}. The latter formula was used in \cite{MTV2} to obtain
a geometric description of the $\B^{K=0}$-action on $(\V^+)^{sing}_\bla$.
Using formula \Ref{grad} we can obtain a similar geometric description of
the $\B^{K=0}$-action on $(\frac 1D \V^-)^{sing}_\bla$.

Let $\bla\in\Z^N_{\geq 0}$, \,$\la_1\geq \dots\geq\la_N, |\bla|=n$.
Introduce $P=\{d_1,\dots,d_N\}$, \,$d_i = \la_i+N-i, i=1,\dots, N$. Let
\be
\Sig_i(u)\,=\,u^{d_i} + \sum_{j=1,\ d_i-j\notin P}^{d_i}\Sig_{ij} u^{d_i-j}.
\ee
Consider the polynomial algebras
\vvn.3>
\be
\C[\bs\Sig]^{sing}:=
\,\C[\Sig_{ij}, i=1,\dots,N, j\in\{1,\dots,d_i\},d_i-j\notin P],
\qquad \C[\bs \si]\,:=\,\C[\si_1,\dots,\si_n].
\vv.2>
\ee
We have
\vvn-.4>
\be
\Wr(\Sig_1(u),\dots,\Sig_N(u))\,=
\prod_{1\le i<j\le N}\!(d_j-d_i)\cdot
\Bigl(u^n+\sum_{s=1}^n\,(-1)^s\>A_s(\bs\Sig)\,u^{n-s}\Bigr)\,,
\vv.3>
\ee
where $A_1(\bs\Sig),\dots,A_n(\bs\Sig)\in\C[\bs\Sig]^{sing}$.
Define an algebra homomorphism
\vvn.3>
\be
\W\,:\,\C[\bs\si]\,\to\,\C[\bs\Sig]^{sing},\qquad
\si_s\,\mapsto\,A_s(\bs\Sig)\,.
\vv.3>
\ee
The homomorphism defines a $\C[\bs\si]$-module structure on
$\C[\bs\Sig]^{sing}$. Define a differential operator $\D_{\bs\Si}$ by
\vvn-.2>
\be
\D_{\bs\Sig}=\,\frac1{\Wr(\Sig_1(u),\dots,\Sig_N(u))}\,\rdet
\left(\begin{matrix} \Sig_1(u) & \Sig_1'(u) &\dots & \Sig_1^{(N)}(u) \\
\Sig_2(u) & \Sig_2'(u) &\dots & \Sig_2^{(N)}(u) \\ \dots & \dots &\dots & \dots \\
1 & \der &\dots & \der^N
\end{matrix}\right).
\vv.2>
\ee
It is a differential operator in the variable $u$, whose coefficients are
formal power series in $u^{-1}$ with coefficients in $\C[\bs\Sig]^{sing}$,
\vvn-.1>
\be
\D_{\bs\Sig}=\,\der^N+\sum_{i=1}^N\,F_i(u)\,\der^{N-i}\>,
\qquad F_i(u)\,=\,\sum_{j=i}^\infty F_{ij}\>u^{-j}\,,
\vv.2>
\ee
and $F_{ij}\in\C[\bs\Sig]^{sing}$, \,$i=1,\dots,N$, \,$j\geq i$.

\begin{thm}
\label{thm 1 0}
The map
\be
\tau_\bla^-\,:\,B^{K=0}_{ij}|_{(\frac 1D\V^-)^{sing}_\bla}\,\mapsto\,F_{ij}
\vv.3>
\ee
defines an isomorphism of the Bethe algebra
$\B^{K=0}\bigl((\frac 1D\V^-)^{sing}_\bla\bigr)$ and the algebra
$\C[\bs\Sig]^{sing}$. The isomorphism $\tau_\bla^-$ becomes an isomorphism
of the $U(\z_N[t])|_{(\frac 1D\V^-)^{sing}_\bla}$-module
$\B^{K=0}\bigl((\frac 1D\V^-)^{sing}_\bla\bigr)$
and the $\C[\bs\si]$-module $\C[\bs\Sig]^{sing}$ if we identify the algebras
$U(\z_N[t])|_{(\frac 1D\V^-)^{sing}_\bla}$ and $\C[\bs\si]$ by the map
$\si_s[\bs z]\mapsto \si_s$, \,$s=1,\dots,n$.
\end{thm}

Fix a vector $v^-\in (\frac 1D\V^-)^{sing}_\bla$ of degree
\;$\sum_{i=1}^N(1-i)\la_i$. By formula \Ref{grad} such a vector is unique up to
proportionality.

\begin{thm}
\label{thm 2 0}
The map
\be
\label{map mu}
\mu_\bla^-\,:\,B^{K=0}_{ij}v^-\,\mapsto\,F_{ij},
\vv.3>
\ee
defines a linear isomorphism $(\frac1D\V^-)^{sing}_\bla\to\C[\bs\Sig]^{sing}$.
The maps $\tau_\bla^-, \mu_\bla^-$ give an isomorphism of the
$\B^{K=0}\bigl((\frac 1D\V^-)^{sing}_\bla\bigr)$-module $(\frac
1D\V^-)^{sing}_\bla$ and the regular representation of the algebra
$\C[\bs\Sig]^{sing}$.
\end{thm}

The proofs of Theorems~\ref{thm 1 0} and~\ref{thm 2 0} are basically
word by word the same as the proofs of Theorems~5.3 and~5.6 in \cite{MTV2}.

\section{Relations with quantum cohomology}
\label{sec quantum}
In lectures \cite{O} Okounkov, in particular, considers the equivariant quantum cohomology
$QH_{GL_n\times \C^*}(T^*F_\bla)$ of the cotangent bundle $T^*F_\bla$ of a flag
variety $F_\bla$.
More precisely, he considers the standard equivariant cohomology
$H^*_{GL_n\times \C^*}(T^*F_\bla)$ as a module
over the algebra of quantum multiplication
and describes this module as the Yangian Bethe algebra of
the {\sl XXX\/} model associated with $V^{\otimes n}$.

The algebra $H^*_{GL_n\times \C^*}(T^*F_\bla)$
has $n+1$ equivariant parameters $z_1,\dots,z_n,u$. The parameters
$z_1,\dots,z_n$ correspond
to the $GL_n$-action on $T^*F_\bla$ and $u$ corresponds of the $\C^*$-action on
$T^*F_\bla$
stretching the cotangent vectors. The operators of quantum multiplication
depend on additional
parameters $q_1,\dots,q_N$ corresponding to quantum deformation.

It is well-known how the Yangian Bethe algebra degenerates into the Gaudin
Bethe algebra, see
for example \cite{T}, \cite{MTV1}. This degeneration
construction gives us the following fact.
Introduce new parameters $K_1,\dots, K_N$ by the formula $q_i=1+K_iu, i=1,\dots,N$, and
consider the limit of the algebra of quantum
multiplication on $H^*_{GL_n\times \C^*}(T^*F_\bla)$ as $u\to 0$. Then this limit is isomorphic to
the $\B^K(\V^+_\bla)$-module $\V^+_\bla$.
This limit is also isomorphic to the $\B^K(\frac 1D\V^-_\bla)$-module $\frac 1D\V^-_\bla$.

\appendix
\section*{Appendix. Topological description of the $\glnt$-module structures\\
on the cohomology of flag varieties}
\refstepcounter{section}

Given $\bla\in \Z_{\geq 0}^N$ define
\begin{align*}
e_{a,a+1}\bla\,&{}=\,
(\la_1,\ldots,\la_{a-1},\la_a+1,\la_{a+1}-1,\la_{a+2},\ldots,\la_N)\,,
\\
e_{a+1,a}\bla\,&{}=\,
(\la_1,\ldots,\la_{a-1},\la_a-1,\la_{a+1}+1,\la_{a+2},\ldots,\la_N)\,,
\\
\bla' &{}=\,
(\la_1,\ldots,\la_{a-1},\la_a,1,\la_{a+1}-1,\la_{a+2},\ldots,\la_N)\,,
\\
\bla'' &{}=\,
(\la_1,\ldots,\la_{a-1},\la_a-1,1,\la_{a+1},\la_{a+2},\ldots,\la_N)\,.
\end{align*}
Let $A'$ (resp. $B', C'$) be the rank $\la_{a}$ (resp. rank $1$, $\la_{a+1}-1$) bundle over $\F_{\bla'}$ whose fiber
over the flag $L_1\subset \ldots \subset L_{N+1}$ is $L_{a}/L_{a-1}$ (resp. $L_{a+1}/L_{a}$, $L_{a+2}/L_{a+1}$).
Let $A''$ (resp. $B'', C''$) be the rank $\la_{a}-1$ (resp. rank $1$, $\la_{a+1}$) bundle over $\F_{\bla''}$ whose fiber
over the flag $L_1\subset \ldots \subset L_{N+1}$ is $L_{a}/L_{a-1}$ (resp. $L_{a+1}/L_{a}$, $L_{a+2}/L_{a+1}$).

Consider the obvious projections
\be
\F_\bla\buildrel{\,\pi'_1\!}\over\longleftarrow\F_{\bla'}
\buildrel{\pi'_2}\over\longrightarrow\F_{e_{a,a+1}\bla}\qquad\text{and}
\qquad\F_\bla\buildrel{\,\pi''_1\!}\over\longleftarrow\F_{\bla''}
\buildrel{\pi''_2}\over\longrightarrow\F_{e_{a+1,a}\bla}\,.
\ee

For an equivariant map $f$ (eg. $f=\pi'_1$ or $\pi''_1$) the induced
pull-back map on equivariant cohomology will be denoted by $f^*$.
For an equivariant fibration $f$ (eg. $f=\pi'_2$ or $\pi''_2$) its Gysin
map (a.k.a.~{\em push-forward} map, or {\em integration along the fibers} map)
will be denoted by $f_*$. The equivariant Euler class of a vector bundle $X$
will be denoted by $e(X)$.

The following theorem was announced in \cite{RSV}.

\begin{thm}
\label{appthm1}
\strut
\begin{enumerate}
\item[(i)]
The map \;$\rho^-(e_{a,a+1}\otimes t^j): H_\bla \to H_{e_{a,a+1}\bla}$
\vvn.2>
\be
x\,\mapsto\,
{\pi'_2}_*\Bigl({\pi'_1}^*(x)\cdot e\bigl(\Hom(B',C')\bigr)\cdot e(B')^j\Bigr)
\vv.2>
\ee
makes the diagram
\vvn-.4>
\be
\begin{CD}
H_\bla @>{\quad\rho^-(e_{a,a+1}\>\otimes\,t^j)\quad}>> H_{e_{a,a+1}\bla}\\
@VV{i^-}V @VV{i^-}V\\[-3pt]
\frac1{D}\V_\bla^- @>{e_{a,a+1}\>\otimes\,t^j}>
\hphantom{{\quad\rho^-(e_{a,a+1}\>\otimes\,t^j)\quad}}>
\frac1{D}\V_{e_{a,a+1}\bla}^-
\end{CD}
\vv-.6>
\ee
commutative.

\medskip
\item[(ii)]
The map \;$\rho^-(e_{a+1,a}\otimes t^j): H_\bla \to H_{e_{a+1,a}\bla}$
\be
x\,\mapsto\,{\pi''_2}_*
\Bigl({\pi''_1}^*(x)\cdot e\bigl(\Hom(A'',B'')\bigr)\cdot e(B'')^j\Bigr)
\ee
makes the diagram
\vvn-.4>
\be
\begin{CD}
H_\bla @>{\quad\rho^-(e_{a+1,a}\>\otimes\,t^j)\quad}>> H_{e_{a+1,a}\bla}\\
@VV{i^-}V @VV{i^-}V\\[-3pt]
\frac1{D}\V_\bla^- @>{e_{a+1,a}\>\otimes\,t^j}>
\hphantom{{\quad\rho^-(e_{a+1,a}\>\otimes\,t^j)\quad}}>
\frac1{D}\V_{e_{a+1,a}\bla}^-
\end{CD}
\vv-.4>
\ee
commutative.
\end{enumerate}
\end{thm}

\begin{proof}
We will prove part (i), the proof of part (ii) is similar. Let $K$ be the index in $\mc I_{e_{a,a+1}\bla}$ with $K_1=\{1,\ldots,(e_{a,a+1}\bla)_1\}$, $K_2=\{(e_{a,a+1}\bla)_1+1,\ldots,(e_{a,a+1}\bla)_1+(e_{a,a+1}\bla)_2\}$, etc.

Consider $x=[h(\bs z, \Gamma_1, \ldots, \Gamma_N)] \in H_\bla$. Its $i^-$-image
is
\vvn.3>
\be
\sum_{I\in\Il}\,v_I \otimes \frac{h(\bs z, \bs z_{I_1}, \ldots, \bs z_{I_N})}
{R(\bs z_{I_1}| \ldots | \bs z_{I_N})}\;.
\ee
The coefficient of $v_K$ of the $e_{a,a+1} \otimes t^j$-image of this is
\vvn.2>
\begin{align}
\label{rr_temp}
&\sum_{i \in K_a}\,
\frac{h(\bs z, \bs z_{K_1}, \ldots, \bs z_{K_{a-1}}, \bs z_{K_a-i},
\bs z_{K_{a+1}\cup i}, \bs z_{K_{a+2}},\ldots, \bs z_{K_N})\, z_i^j}
{R( \bs z_{K_1}, \ldots, \bs z_{K_{a-1}}, \bs z_{K_a-i}, \bs z_{K_{a+1}\cup i},
\bs z_{K_{a+2}},\ldots, \bs z_{K_N} )}\,={}
\\[3pt]
&\,\,{}=\,\frac1{R(\bs z_{K_1}|\ldots | \bs z_{K_N})}\,\sum_{i \in K_a}\,
\frac{h(\bs z, \bs z_{K_1}, \ldots, \bs z_{K_a-i}, \bs z_{K_{a+1}\cup i},
\ldots, \bs z_{K_N})\,z_i^j\,R(z_i| \bs z_{K_{a+1}})}
{R( \bs z_{K_a-i}, z_i ) }\,.\!\!
\notag
\end{align}
On the other hand, the $\rho^-(e_{a,a+1}\otimes t^j)$-image of $x$ (using a version of the Atiyah-Bott localization formula for ${\pi'_2}_*$) is
\vvn.3>
\beq
\label{rr}
\sum_{\delta\in \Delta_a}
\frac{h(\bs z, \Delta_1, \ldots, \Delta_{a-1},\Delta_a-\delta, \delta, \Delta_{a+1},\ldots,\Delta_N) R(\delta|\Delta_{a+1})\delta^j }
{R(\Delta_a-\delta | \delta)},
\vv.1>
\eeq
where we denoted the Chern roots of the natural bundles over
$\F_{e_{a,a+1}\bla}$ by $\Delta_1,\ldots, \Delta_N$.
The coefficient of $v_K$ of
%its
$i^-$-image of \Ref{rr} is the right hand side of \Ref{rr_temp}.
The theorem is proved.
\end{proof}

The topological interpretation of generators of the $\rho^+$-representation is
similar, its proof is left to the reader.

\goodbreak
\begin{thm}
\label{appthm2}
\strut
\begin{enumerate}
\item[(i)]
For the map \;$\rho^+(e_{a,a+1}\otimes t^j):H_\bla\to H_{e_{a,a+1}\bla}$,
\vvn.2>
\be
x\,\mapsto\,
{\pi'_2}_*\Bigl({\pi'_1}^*(x)\cdot e\bigl(\Hom(A',B')\bigr)\cdot e(B')^j\Bigr)
\vv.2>
\ee
we have \;\,$i^+\!\circ\rho^+(e_{a,a+1}\otimes t^j)\>=\>
(e_{a,a+1}\otimes t^j)\circ i^+$.

\medskip
\item[(ii)]
For the map \;$\rho^+(e_{a+1,a}\otimes t^j):H_\bla\to H_{e_{a+1,a}\bla}$,
\vvn.2>
\be
x \mapsto{\pi''_2}_*\left({\pi''_1}^*(x) \cdot e(\Hom(B'',C''))\cdot e(B'')^j \right)
\vv.2>
\ee
we have \,\;$i^+\!\circ\rho^+(e_{a+1,a}\otimes t^j)\>=
\>(e_{a+1,a}\otimes t^j)\circ i^+$.
\end{enumerate}
\end{thm}

\medskip
The $\glnt$-module structures $\rho^\pm$
on $\bigoplus_\bla H_\bla$ descend to $\glnt$-module structures on $H(\C)$,
also denoted by $\rho^\pm$
in Section \ref{rr_sec}. The topological interpretation of the actions of
$e_{a,a+1}\otimes t^j$ and $e_{a+1,a}\otimes t^j$ for these representations
is the same as that
for $\bigoplus_\bla H_\bla$ given in Theorems~\ref{appthm1} and~\ref{appthm2}.

\medskip

Topological description of the $\glnt$-module structures on the cohomology of
flag varieties, presented in this Appendix, is a version of a construction due
to Ginzburg and Vasserot, cf.~\cite{Vas1}, \cite{Vas2}. We thank Eric Vasserot
for pointing out to us these references. However, unlike these authors, we work
with flag varieties themselves and not with their cotangent bundles. Also our
proof is different.


\bigskip
\begin{thebibliography}{[CG]}
\normalsize
\frenchspacing
\raggedbottom

\bi[AB]{AB} M.\,F.\,Atiyah, R.\,Bott, {\it
The moment map and equivariant cohomology}, Topology {\bf 23} (1984), 1--28

\bi[BMO]{BMO}
A.\,Braverman, D.\,Maulik, A.\,Okounkov, {\it
Quantum cohomology of the Springer resolution}, Preprint (2010), 1--35, {\tt
arXiv:1001.0056}

\bi [CL]{CL} V.\,Chari, S.\,Loktev, {\it Weyl, Fusion and Demazure
modules for the current algebra of $\frak{sl}_{r+1\/}$}, Adv. Math. 207
(2006), no.\;2, 928--960

\bi [CP]{CP} V.\,Chari, A.\,Pressley
{\it Weyl Modules for Classical and Quantum Affine algebras\/},
Represent. Theory 5 (2001), 191--223 (electronic)

\bi [CT]{CT} A.\,Chervov, D.\,Talalaev,
{\it Quantum spectral curves, quantum integrable systems and the geometric
Langlands correspondence\/}, Preprint (2006), 1--54; {\tt hep-th/0604128}

\bi[G1]{G1} M.\,Gaudin, {\it Diagonalisation d'une classe d'Hamiltoniens de spin}, J. Physique {\bf 37} (1976), no. 10, 1089--1098

\bi[G2]{G2} M.\,Gaudin, {\it La fonction d'onde de Bethe}, Collection du
Commissariat \`a l'\'Energie Atomique: S\'erie Scientifique, Masson, Paris,
1983

\bi [MTV1]{MTV1}
E.\,Mukhin, V.\,Tarasov, A.\,Varchenko,
{\it Bethe Eigenvectors of Higher Transfer Matrices\/},
J.~Stat. Mech. (2006), no.\;8, P08002, 1--44

\bi [MTV2]{MTV2}
E.\,Mukhin, V.\,Tarasov, A.\,Varchenko, {\it
Schubert calculus and representations of general linear group},
J.\;Amer. Math. Soc.~{\bf 22} (2009), no.\;4, 909--940

\bi [MTV3]{MTV3} E.\,Mukhin, V.\,Tarasov, A.\,Varchenko, {\it
Spaces of quasi-exponentials and representations of $\gln$},
J.\;Phys.~A {\bf 41} (2008), no.\;19, 194017, 28 pp.

\bibitem[MV]{MV}
E.\,Mukhin, A.\,Varchenko,
{\it Critical points of master functions and flag varieties\/},
Commun. Contemp. Math. {\bf 6} (2004), no.\;1, 111--163

\bi[O]{O} A.\,Okounkov, {\it Quantum Groups and Quantum Cohomology}, Lectures
at the 15th Midrasha Mathematicae on: "Derived categories of algebro-geometric
origin and integrable systems", December 19--24, 2010, Jerusalem

\bi[RV]{RV} R.\,Rim\'anyi, A.\,Varchenko, {\it Conformal blocks in the tensor
product of vector representations and localization formulas}, Preprint (2009),
1--21, {\tt arXiv:0911.3253}

\bi[RSV]{RSV} R.\,Rim\'anyi, V.\,Schechtman, A.\,Varchenko,
{\it Conformal blocks and equivariant cohomology}, Preprint (2010), 1--23,
{\tt arXiv:1007.3155}

\bi [T]{T} D.\,Talalaev,
{\it Quantization of the Gaudin System\/}, Preprint (2004), 1--19;\\
{\tt hep-th/0404153}

\bi[V]{V} A.\,Varchenko, {\it
A Selberg integral type formula for an $\,\frak{sl}_2$ one-dimensional space
of con\-formal blocks}, Mosc. Math.~J. {\bf 10} (2010), no. 2, 469--475, 480

\bi[Vas1]{Vas1} E.\,Vasserot, {\it Repr\'esentations de groupes quantiques et
permutations\/}, Ann. Sci. ENS {\bf 26} (1993), 747--773

\bi[Vas2]{Vas2} E.\,Vasserot, {\it Affine quantum groups and equivariant
$K$-theory}, Transformation groups, {\bf 3} (1998), 269--299

\end{thebibliography}
\end{document}